\documentclass[12pt]{article}
\usepackage{amsmath, amssymb, theorem}

\theoremheaderfont\scshape
\newtheorem{theorem}{Theorem}[section]

\newtheorem{corollary}[theorem]{Corollary}
\newtheorem{definition}[theorem]{Definition}
\newtheorem{open}[theorem]{Problem}
\newtheorem{proposition}[theorem]{Proposition}

\theorembodyfont{\rmfamily}

\newcommand\comment[1]{}

\newcommand{\nohyphens}{\hyphenpenalty=10000\exhyphenpenalty=10000\relax}%
\makeatletter
\renewcommand{\@seccntformat}[1]{\csname the#1\endcsname.\hspace{1em}}%
\renewcommand\section{\@startsection{section}{1}{\z@}%
                {-3.5ex \@plus -1ex \@minus -.2ex}%
                {2.3ex \@plus.2ex}%
                {\setcounter{equation}0\bfseries\nohyphens}}%
\makeatother

\begin{document}

\title{A Short Note on some open problems in the geometry of operator ideals}
\date{}
\maketitle

\section{The principle of local reflexivity for maximal Banach ideals}
Throughout this Short Note, we essentially adopt notation and
terminology from \cite{p1}, which however cannot be explained here
in detail, due to the limitation of space. Therefore, we would like
to refer the interested reader to \cite{df, j, p1}, the excellent
survey article \cite{djp} and the further references therein.

Let $(\frak{A},\mathbf{A})$ be an arbitrary \textit{maximal} Banach
ideal. It seems to be \textit{still} an open problem whether it is
\textit{always} possible to transfer the norm estimation in the
far-reaching classical principle of local reflexivity to the ideal
norm $\mathbf{A}$. If this were not the case, we would be very
interested in constructing an explicit counterexample. More
precisely, we would like to know whether there exists a
\textit{maximal} Banach ideal $(\frak{A}_0,\mathbf{A}_0)$ which does
\textit{not} satisfy the following factorization property:
\begin{definition}[cf. \cite{oe1, oe2}]
Let $E$ and $Y$ be Banach spaces, where $E$ is finite dimensional
and $F$ is a finite dimensional subspace of $Y^\prime$. Let
$(\frak{A}, \mathbf{A})$ be a maximal Banach ideal and $\varepsilon
>0$. We say that the principle of $\frak{A}$-local reflexivity
$($short: $\frak{A}$-LRP$)$ is satisfied, if for every $T\in
\frak{L}(E,Y^{\prime \prime })$ there exists an operator $S \in
\frak{L}(E,Y)$ such that
\begin{enumerate}
\item[$(i)$]  $\mathbf{A}(S)\leq (1+\epsilon )\cdot \mathbf{A}(T)$;
\item[$(ii)$]  $\left\langle Sx,y^{\prime }\right\rangle =\left\langle
y^{\prime },Tx\right\rangle $ for all $(x,y^{\prime })\in E\times
F$;
\item[$(iii)$]  $j_YSx=Tx$ for all $x\in T^{-1}(j_Y(Y))$.\footnote{Here, $j_Y
: Y \hookrightarrow Y^{\prime\prime}$ denotes the canonical
embedding from the Banach space $Y$ into its bidual
$Y^{\prime\prime}$.}
\end{enumerate}
\end{definition}
Obviously, the classical principle of local reflexivity is simply
the $\frak{B}$-LRP, where $(\frak{B}, \mathbf{B}) : = (\frak{L},
\Vert \cdot \Vert)$.
\begin{open}\label{OP 1}
Does every maximal Banach ideal satisfy the $\frak{A}$-LRP?
\end{open}
Due to the finite nature of \textit{maximal} Banach ideals and the
$\frak{A}$-LRP, local versions of injectivity and surjectivity of
operator ideals play a key role. Moreover, they imply interesting
results for operators with infinite dimensional range (cf.
\cite{oe3, oe4, oe5}). These are the following so-called
``accessibility conditions'', treated in detail in \cite{de, df}:
\begin{definition}
Let $(\frak{B},\mathbf{B})$ be a $p$-normed Banach ideal, where $0
<p \leq 1$.
\begin{itemize}
\item[$($i$)$] $(\frak{B},\mathbf{B})$ is called right-accessible, if for
any two Banach spaces $E$ and $Y$, $\mbox{dim}(E) < \infty$, any
operator $T\in \frak{L}(E,Y)$ and any $\varepsilon >0$ there are a
finite dimensional subspace $N$ of $Y$ and $S\in \frak{L}(E,N)$ such
that $T=J_N^YS$ and $\mathbf{B}(S)\leq (1+\varepsilon
)\mathbf{B}(T)$; here $J_N^Y : N \hookrightarrow Y$ denotes the
canonical embedding.
\item[$($ii$)$] $(\frak{B},\mathbf{B})$ is called left-accessible, if for
any two Banach spaces $F$ and $X$, $\mbox{dim}(F) < \infty$, any
operator $T\in \frak{L}(X,F)$ and any $\varepsilon >0$ there are a
finite codimensional subspace $L$ of $X$ and $S\in \frak{L}(X/L,F)$
such that $T=SQ_L^X$ and $\mathbf{B}(S)\leq (1+\varepsilon
)\mathbf{B}(T)$; here $Q_L^X : X \twoheadrightarrow X/L$ denotes the
canonical quotient map.
\item[$($iii$)$] $(\frak{B},\mathbf{B})$ is called
accessible, if it is both, right-accessible and left-accessible.
\end{itemize}
\end{definition}
\begin{proposition}[cf. \cite{oe5}]\label{right-acc implies A-LRP}
If $(\frak{A},\mathbf{A})$ is a right-accessible maximal Banach
ideal, then the $\frak{A}$-LRP holds.
\end{proposition}
We \textit{still} do not know whether the statement in Proposition
\ref{right-acc implies A-LRP} can be reversed. In 1993, Pisier
constructed explicitly a maximal Banach ideal $(\frak{A}_P,
\mathbf{A}_P)$ which is not right-accessible, and solved - 37 years
later only - a further problem of Grothendieck's famous
R\'{e}sum\'{e} (cf. \cite{gr} and \cite{df}). Consequently, Problem
\ref{OP 1} is even harder than Grothendieck's tough accessibility
problem. A maximal Banach ideal $(\frak{A},\mathbf{A})$ which does
not satisfy the $\frak{A}$-LRP, necessarily cannot be
right-accessible! Notice that in the investigation of those
problems, Banach spaces \textit{without} the approximation property
(such as the Pisier space $P$) are necessarily involved (cf.
\cite{df, oe1, oe5}). Consequently, we arrive at
\begin{open}\label{OP 2}
Let $(\frak{A},\mathbf{A})$ be a maximal Banach ideal. Are then the
following statements equivalent?
\begin{itemize}
\item[$($i$)$] $(\frak{A},\mathbf{A})$ is right-accessible.
\item[$($ii$)$] $(\frak{A},\mathbf{A})$ satisfies the $\frak{A}$-LRP.
\end{itemize}
\end{open}
\begin{open}\label{OP 3}
Does $(\frak{A}_P, \mathbf{A}_P)$ satisfy the $\frak{A}_P$-LRP?
\end{open}
Let $(\frak{B},\mathbf{B})$ be a $p$-normed Banach ideal, where $0
<p \leq 1$. In 1971, Pietsch defined the adjoint operator ideal
$(\frak{B}^\ast, \mathbf{B}^\ast )$ to uncover the structure of
maxi-mal Banach ideals.
In 1973, Gordon-Lewis-Retherford introduced a related smaller
``conjugate'' operator ideal $(\frak{B}^\Delta, \mathbf{B}^\Delta
)$, which - in the contrary - is somehow ``skew'' under the point of
view of trace duality but nevertheless quite useful (cf. \cite{glr,
oe1}):
\begin{definition}
Let $X$, $Y$ two Banach spaces and $\frak{B}^\Delta (X,Y)$ be the
set of all those operators $T\in \frak{L}(X,Y)$ satisfying
\[
\mathbf{B}^\Delta (T):=\sup \big\{\vert \textup{tr}(TL) \vert : L\in
\frak{F}(Y,X),\mathbf{B}(L)\leq 1\big\}< \infty .
\]
Then $(\frak{B}^\Delta ,\mathbf{B}^\Delta )$ is a Banach
ideal.\footnote{Here $\mbox{tr} : \frak{F}(Y,Y) \longrightarrow
{\Bbb{K}}$ denotes the trace on the operator ideal component of all
finite rank operators on $Y$.} It is called the \textit{conjugate
ideal} of $(\frak{B},\mathbf{B})$.
\end{definition}
It is easy to see that the various approximation properties of
Banach spaces and the accessibility of operator ideals are
intrinsically related to a ``good behaviour'' of trace duality. If
$\frak{A}$ is a maximal Banach ideal, then $(\frak{A}^\Delta
,\mathbf{A}^\Delta )$ is right-accessible (cf. \cite{oe1, oe2}).
However, we do not know whether $\frak{A}^\Delta$ is
left-accessible, since:
\begin{theorem}[cf. \cite{oe1, oe2}]
Let $(\frak{A},\mathbf{A})$ be a maximal Banach ideal. Then the
following statements are equivalent:
\begin{enumerate}
\item[$(i)$]  $\frak{A}^\Delta $ is left--accessible.
\item[$(ii)$]  $\frak{A}(E,Y^{\prime \prime })$ $\widetilde{=}$ $\frak{A}(E,Y)^{\prime \prime }$
for all finite dimensional Banach spaces $E$ and arbitrary Banach
spaces $Y$.
\item[$(iii)$]  The $\frak{A}$-LRP holds.
\end{enumerate}
\end{theorem}
\section{Property (I) and normed products of Banach ideals}
After extending finite rank operators in suitable quasi-Banach
ideals $\frak{B} \not= \frak{L}$, the $\frak{B}$-LRP and the
calculation of conjugate ideal-norms allow us to ignore the
structure of the range space and to leave the finite dimensional
case. There are sufficient conditions on $\frak{B}$ which imply that
for any two Banach spaces $X$ and $Y$, any Banach space $Z$ which
contains $X$ isometrically, any finite rank operator $L \in
\frak{F}(X, Y)$ and any $\varepsilon > 0$ there exists a
\textit{finite rank extension} $\widetilde{L} \in \frak{F}(Z, Y)$
such that $\mathbf{B}(\widetilde{L}) \leq (1+\epsilon )\cdot
\mathbf{B}(L)$ (cf. \cite{oe5}).\footnote{Recall that there is no
Hahn-Banach extension theorem for finite rank operators, viewed as
elements of the Banach ideal $(\frak{L},\Vert \cdot \Vert)$.}

Consequently, we have to look for a suitable class of product ideals
$\frak{A}\circ \frak{B}$ which satisfy the following factorization
property: Given $\varepsilon
> 0$, any finite rank operator $L \in \frak{A}\circ \frak{B}$ can be
factorized as $L=AB$ such that $\mathbf{A}(A)\cdot \mathbf{B}(B)\leq
(1+\varepsilon )\cdot \mathbf{A} \circ \mathbf{B}(L)$ and $A$
respectively $B$ has finite dimensional range. Product ideals
$\frak{A}\circ \frak{B}$ of this type are said to have
\textit{property $(I)$} respectively \textit{property $(S)$}. They
had been introduced in \cite{jo} to prepare a detailed investigation
of trace ideals.

In view of looking for a counterexample of a maximal Banach ideal $(%
\frak{A}_0,\mathbf{A}_0)$ which does not satisfy the
$\frak{A}_0$-LRP, property (I) of product ideals of type
$\frak{A}^{*}\circ \frak{L}_\infty $ seemingly plays a key role.
\begin{theorem}[cf. \cite{oe5}]\label{injective hulls in L2}
Let $(\frak{A},\mathbf{A})$ be a maximal Banach ideal such that
$\frak{A}^{*}\circ \frak{L}_\infty $ has property $(I)$. If
\textup{space}$(\frak{A})$ contains a
Banach space $X_0$ such that $X_0$ has the bounded approximation property but $%
X_0^{\prime \prime }$ has not, then the $\frak{A}^{*}$-LRP is not
satisfied.
\end{theorem}
A further interesting class is given by the family of all operator
ideals which contain $(\frak{L}_2,\mathbf{L}_2)$ as a factor
(i.\,e., the maximal and injective Banach ideal of all operators
which factor through a Hilbert space). This class does not only show
surprising connections with the principle of local reflexivity for
operator ideals. There are also links to the existence of an
ideal-\textit{norm} on certain product ideals and the property (I),
reflected e.\,g. in the following two results (cf. \cite{oe6}):
\begin{theorem}
Let $(\frak{A},\mathbf{A})$ be an arbitrary maximal Banach ideal
such that $\frak{A} \circ \frak{L}_2$ is normed. If the $(\frak{A}
\circ \frak{L}_2)^{\ast\ast}$-LRP is satisfied, then $(\frak{A}
\circ \frak{L}_2)^{\ast\ast} \circ \frak{L}_\infty$ has both,
property $(I)$ and property $(S)$.
\end{theorem}
\begin{theorem}\label{normed L2 factor and LRP}
Let $(\frak{A},\mathbf{A})$ be an arbitrary maximal Banach ideal
such that $\frak{A} \circ \frak{L}_2$ is normed. If $\frak{L}_2
\subseteq \frak{A}^\ast$, then the $(\frak{A} \circ
\frak{L}_2)^{\ast\ast}$-LRP cannot be satisfied.
\end{theorem}
Here is a nice application of Theorem \ref{normed L2 factor and
LRP}: Put $\frak{A} : = (\frak{L}_1 \circ
\frak{L}_\infty)^{\ast\ast}$. Due to Grothendieck's inequality in
operator form, it follows that $\frak{L}_2 \subseteq \frak{A}^\ast$
(cf. \cite{df, oe4}). Hence, in view of Theorem \ref{normed L2
factor and LRP}, a natural question appears:
\begin{open} Is the product ideal $(\frak{L}_1 \circ \frak{L}_\infty)^{\ast\ast} \circ
\frak{L}_2$ normed?
\end{open}
However, we do not know criteria which are sufficient for the
existence of an ideal-\textit{norm} on a given product of
quasi-Banach ideals. Here, one has to study very carefully the
structure of $F$-spaces which are not locally convex, such as
$L^p([0,1])$ and the Hardy space $H^p$ on the unit disk, for some
$0<p<1$ (cf. \cite{kpr}). It seems to be much easier to provide
arguments which imply the non-existence of such an ideal-norm (by
using trace ideals) (cf. \cite{oe6}).

\begin{theorem}[cf. \cite{oe4, oe5}]
Let $(\frak{A},\mathbf{A})$ be a maximal Banach ideal such that
$\frak{A}^{*}\circ \frak{L}_\infty $ is right-accessible and has
property $(I)$. Let $X$ and $Y$ be Banach spaces such that
$X^\prime$ and $Y$ are of cotype 2. If the $\frak{A}^\ast$-LRP is
satisfied, then
\[
\frak{A}^{\rm{\textup{inj}}}(X,Y)\subseteq \frak{L}_2(X,Y)
\]
and
\[
\mathbf{L}_2(T)\leq \big(2\mathbf{C}_2(X^{\prime })\cdot
\mathbf{C}_2(Y)\big)^{\frac{3}{2}}\cdot
\mathbf{A}^{\rm{\textup{inj}}}(T)
\]
for all operators $T\in \frak{A}^{\rm{\textup{inj}}}(X,Y)$.
\end{theorem}
\begin{corollary}\label{right-acc and LRP and prop (I)}
Let $(\frak{B},\mathbf{B})$ be a maximal Banach ideal. Let $X_0$ and
$Y_0$ be Banach spaces such that $X_0^\prime$ and $Y_0$ are of
cotype 2 and $\frak{B}^\ast(X_0,Y_0) \not\subseteq
\frak{L}_2(X_0,Y_0)$. If $\frak{B}\circ \frak{L}_\infty$ is
right-accessible and has property $(I)$, then the $\frak{B}$-LRP is
not satisfied.
\end{corollary}
Consequently, an application of Proposition \ref{right-acc implies
A-LRP} immediately implies the following surprising result:
\begin{corollary}\label{right-acc vs prop (I)}
Let $(\frak{B},\mathbf{B})$ be a maximal Banach ideal. Let $X_0$ and
$Y_0$ be Banach spaces such that $X_0^\prime$ and $Y_0$ are of
cotype 2 and $\frak{B}^\ast(X_0,Y_0) \not\subseteq
\frak{L}_2(X_0,Y_0)$. If $\frak{B}$ is right-accessible, then
$\frak{B}\circ \frak{L} _\infty$ does not have property $(I)$.
\end{corollary}
\begin{open}\label{OP 4}
Let $(\frak{C}_2,\mathbf{C}_2)$ denote the maximal injective Banach
ideal of all cotype 2 operators. Does $\frak{C}_2^{*}\circ
\frak{L}_\infty $ have property $(I)$?
\end{open}
If Problem \ref{OP 4} had a positive answer, Corollary
\ref{right-acc vs prop (I)} would imply that
$(\frak{C}_2,\mathbf{C}_2)$ cannot be left-accessible; and a further
open problem of Defant and Floret could be solved (see \cite{df},
21.2., p. 277).
\begin{open}
Can we drop the assumption ``\,$\frak{B}\circ \frak{L}%
_\infty$ is right-accessible'' in Corollary \ref{right-acc and LRP
and prop (I)}?
\end{open}

\end{document}